\documentclass[leqno,12pt]{amsart}
\usepackage{amsfonts}
\usepackage{amsmath,amssymb,amsthm}
\usepackage{amsfonts}

\newtheorem{theorem}{Theorem}[section]
\newtheorem{proposition}[theorem]{Proposition}
\newtheorem{lemma}[theorem]{Lemma}
\newtheorem{corry}[theorem]{Corollary}

\theoremstyle{definition}

\theoremstyle{definition}
\newtheorem{rem}[theorem]{Remark}

\numberwithin{claim}{theorem}

\renewenvironment{proof}{\textit{Proof.}}{\hfill\ensuremath{\qed}}

\def \qed{\hfill{\hbox{$\square$}}}

\numberwithin{equation}{section}

\begin{document}
\title[Surfaces in $L^4_1(f,c)$]{Surfaces in Robertson-Walker Space-Times with Positive Relative Nullity}

\author[B. Bekta\c{s} Dem\.{i}rc\.i]{Burcu Bekta\c s Dem\.{i}rc\.i}
\address{Department of Software Engineering, Faculty of Engineering, Topkap{\i} Campus, Fatih Sultan Mehmet Vak{\i}f University, 
, \.{I}stanbul, T{\"u}rk\.{I}ye}
\email{bbektas@fsm.edu.tr, 0000-0002-5611-5478}

\author[N. Cenk Turgay]{Nurettin Cenk Turgay}
\address{Department of Mathematics, Faculty of Science and Letters, Istanbul Technical University, \.{I}stanbul, T{\"u}rk\.{I}ye}
\email{nturgay@itu.edu.tr, 0000-0002-0171-3876}

\begin{abstract}
In this article, we study space-like and time-like surfaces in a Robertson-Walker space-time,, denoted by $L^4_1(f,c)$,
having positive relative nullity. 
First, we give the necessary and sufficient conditions for such space-like and time-like surfaces
in $L^4_1(f,c)$. 
Then, we obtain the local classification theorems for space-like and time-like surfaces in $L^4_1(f,0)$ with 
positive relative nullity.
Finally, we consider the space-like and time-like surfaces 
in $\mathbb{E}^1_1\times\mathbb{S}^3$ and $\mathbb{E}^1_1\times\mathbb{H}^3$ with positive relative nullity.
These are the special spaces of $L^4_1(f,c)$ when the warping function $f$ is a constant function, with $c=1$ 
for $\mathbb{E}^1_1\times\mathbb{S}^3$ and $c=-1$ for $\mathbb{E}^1_1\times\mathbb{H}^3$.
\end{abstract}


\maketitle
\section{Introduction}
The relative null space of a submanifold of a pseudo-Riemannian manifold at a point is the subspace of its tangent space on which the second fundamental form vanishes identically.
Then, a submanifold in a pseudo-Riemannian manifold
is said to have positive relative nullity
if, at every point on the submanifold, the dimension of the relative null space is greater than zero.

There are several works related to submanifolds with the relative null space and positive relative nullity
in different ambient spaces. 
M. Dajczer and D. Gromoll introduced the concept of submanifolds with positive relative nullity and 
they obtained the necessary and sufficent conditions for a spherical submanifold to
have positive relative nullity in \cite{DajczerGro}. 
B.-Y. Chen and J. Van der Veken provided a complete classification of marginally trapped surfaces
with positive relative nullity in Lorentzian space-forms in \cite{ChenJVanderVeken2006}. 
M. Dajczer and R. Tojeiro gave some results about submanifolds of space forms having positive relative nullity, 
\cite{DajczerTojeiro2019}.
Additionally, the authors \cite{DemirciTurgay2020} studied marginally trapped surfaces with positive relative nullity 
in 4-dimensional pseudo--Riemannian space forms with index 2.  
 
Another example of the Lorentzian space is the family of Robertson-Walker space times which is an
important family of cosmological models in general relativity.
Some studies have focused on classifying surfaces or hypersurfaces of Robertson-Walker spacetimes
having different geometrical properties such as \cite{AnciauxCipriani2020}, \cite{ChenJVanderVeken2007}, 
\cite{DekimpeJVanderVeken2020}. 
Especially, within the context of positive relative nullity 
B.-Y. Chen and J. Van der Veken \cite{ChenJVanderVeken2007} investigated a marginally trapped surface with positive relative nullity 
in a Robertson Walker space-time and 
they proved that a Robertson Walker space-time, 
which contains no open subsets of constant curvature, does not allow for 
the existence of a marginally trapped surface with positive relative nullity. 

The aim of this paper is to examine space-like and time-like surfaces in a Robertson-Walker spacetime
$L^4_1(f,c)$ that possess positive relative nullity, without assuming the mar\-gi\-nal\-ly trapped condition.
First, we give the necessary and sufficient conditions for such space-like and time-like surfaces
in $L^4_1(f,c)$ and we also get the forms of Levi-Civita connections of them. 
Then, we obtain the classification theorems for space-like and time-like surfaces in $L^4_1(f,0)$ having 
positive relative nullity.
We end the paper with the results concerning space-like and time-like surfaces 
in $\mathbb{E}^1_1\times\mathbb{S}^3$ and $\mathbb{E}^1_1\times\mathbb{H}^3$.


\section {Preliminaries}
In this section, we give a brief summary of basic facts and equations of the theory of submanifolds and Robertson Walker Spaces, 
(see \cite{OneillBook}, \cite{ChenBook}).

\subsection{Robertson-Walker Spaces}
Let $R^{n-1}(c)$ denote the $n-1$ dimensional Riemannian space-form with the constant sectional curvature $c$, i.e.,
$$R^{n-1}(c)=\left\{\begin{array}{cc}
\mathbb S^{n-1}&\mbox{if }c=1,\\
\mathbb E^{n-1}&\mbox{if }c=0,\\
\mathbb H^{n-1}&\mbox{if }c=-1
\end{array}\right.$$
and $g_c$ stand for its metric tensor and $\nabla^{R^3(c)}$ denote its Levi-Civita connection. 
If $I$ is an open interval and $f:I\to\mathbb R$ is a smooth, non-vanishing function, 
then the Robertson-Walker space $L^n_1(f,c)$ is defined 
as the Lorentzian warped product $I^1_1\times_f R^{n-1}(c)$ whose metric tensor $\widetilde g$ is
$$\widetilde g=\langle\cdot,\cdot\rangle=-dt^2+f(t)^2g_c.$$
If the warping function $f$ is constant, then $L^n_1(f,c)$ 
is the Lorentzian product of $I^1_1$ and a Riemannian real space form. 
For any function $f$ and $n = 4$, the Lorentzian manifolds $L^4_1(f,c)$ are known as Robertson-Walker space-times.
Also, the vector field $\displaystyle\frac{\partial}{\partial t} $ is called the comoving observer field in general relativity, \cite{ChenJVanderVeken2007}.

Let  $\Pi^1:I\times R^{n-1}(c)\to I$ and $\Pi^2:I\times R^{n-1}(c)\to R^{n-1}(c)$
denote  the canonical projections. For a given vector field $X$ in ${L}^n_1(f,c)$, we define a function $X_0$ and a vector field $\bar X$ by the orthogonal decomposition
$$X=X_0\frac{\partial}{\partial t} +\bar X,$$
that is,
$$X=X_0\frac{\partial}{\partial t} +\sum\limits_{i=1}^{n-1}X_i\frac{\partial}{\partial x_i}:=(X_0,X_1,\hdots,X_{n-1}),$$ where 
$(x_1,x_2,\hdots,x_{n-1})$ is a Cartesian coordinate system in $\mathbb R^{n-1}$ and we have 
$$X_0=-\langle\frac{\partial}{\partial t} ,X\rangle,\qquad \Pi_1^*(\bar X)=0.$$ First, we would like to express the Levi-Civita connection of $ L^n_1(f,c)$. 
Occasionally, by misusing terminology, we are going to put $(0,\bar X)=\bar X$.

We are going to use the following lemma which can be directly obtained from  \cite{OneillBook} (See also \cite[Lemma 2.1]{ChenJVanderVeken2007}).
\begin{lemma}\label{LemmaLn1f0LCConnect}
The Levi-Civita connection $\widetilde\nabla$ of ${L}^n_1(f,c)$ is
\begin{eqnarray}
\label{RWSTconn}
\widetilde{\nabla}_X Y&=&\nabla^0_XY+\frac{f'}{f} \left(g(\bar X,\bar Y){\frac{\partial}{\partial t} }+ X_0\bar Y+Y_0\bar X\right)
\end{eqnarray}
whenever $X$ and $Y$ are tangent to ${L}^n_1(f,c)$, where  $\nabla^0$ denotes the Levi-Civita connection of the Cartesian product space ${L}^n_1(1,c)=I\times R^{n-1}(c)$.
\end{lemma}

On the other hand, the Riemann curvature tensor $\widetilde R$ of ${L}^n_1(f,c)$ is as follows:
\begin{lemma}\label{LemmaLn1f0CurvTensor}\cite{ChenJVanderVeken2007}.
The Riemannian curvature tensor $\widetilde R$ of ${L}^n_1(f,c)$ satisfies
\begin{align}
\begin{split}
\widetilde R (\frac{\partial}{\partial t} ,\bar X )\frac{\partial}{\partial t} =\frac{f''}f\bar X, \qquad& \widetilde R (\frac{\partial}{\partial t} ,\bar X )\bar Y=\frac{f''}f\langle \bar X,\bar Y\rangle\frac{\partial}{\partial t} ,\\
\widetilde R (\bar X,\bar Y )\frac{\partial}{\partial t} =0, \qquad& \widetilde R (\bar X,\bar Y )\bar Z=\frac{f'{}^2+c}{f^2} (\langle \bar Y,\bar Z\rangle \bar X-\langle \bar X,\bar Z\rangle \bar Y )
\end{split}
\end{align}
whenever $\Pi_1^*(\bar X)=\Pi_1^*(\bar Y)=\Pi_1^*(\bar Z)=0$
\end{lemma}

\subsection{Surfaces in Robertson-Walker Space-Times}
Let $M$ be a non-degenerated surface in $L^4_1(f,c)$, $\nabla$ and $g$ denote its Levi-Civita connection and metric tensor, respectively. One can define the second fundamental form  $ h $, the shape operator  $ A $, and the normal connection  $ \nabla^\perp $  of  $ M $ by the Gauss and Weingarten formulas  given by
\begin{eqnarray}
    \label{Gauss} \widetilde\nabla_X Y &=& \nabla_X Y + h(X,Y), \\
    \label{Weingarten} \widetilde\nabla_X \xi &=& -A_\xi(X) + \nabla^\perp_X \xi,
\end{eqnarray}
where  $ X $  and  $ Y $  are tangent vectors to  $ M $  and  $ \xi $  is a normal vector to  $ M $. Note that the equation
\begin{eqnarray}
    \label{AhRelatedBy} 
    \langle h(X,Y),\xi\rangle = \langle A_\xi X,Y\rangle
\end{eqnarray}
is satisfied.

On the other hand, one can define a tangent vector field $T$ and a normal vector field $\eta$ on $M$ by the decomposition
\begin{equation}
    \label{RWTS-etaandTDef} 
    \left. \frac{\partial}{\partial t}  \right|_M = T + \eta,
\end{equation}
(see \cite{MendonTojeiro2014}, \cite{Tojeiro2010}).

The relative null space $\mathcal{N}_p$ at a $p\in M$ is defined by 
$$\mathcal{N}_p=\{X_p\in T_pM\;|\;h_p(X_p,Y_p)=0 
\mbox{ for all $Y_p\in T_pM$ }\}.$$
If $\dim\mathcal{N}_p>0$ for all $p\in M$, then
$M$ is said to have positive relative nullity, \cite{ChenJVanderVeken2007}.

\begin{rem} 
If $T = 0 $  on a non-empty connected open subset $\mathcal{O} $  of $M$, then
$ \mathcal{O} $  is a horizontal slice, i.e., $  (\mathcal{O}, g) \subset \{ t_0 \} \times_{f(t_0)} {R}^3(c)$. In this case, by considering \cite[Lemma 3.2]{ChenJVanderVeken2007}, one can conclude that $ (\mathcal{O}, g) $ has positive relative nullity if and only if $f'(t_0)=0$ and $\mathcal{O}=\{ t_0 \} \times_{f(t_0)}\tilde M$, where $\tilde M$ is a surface in ${R}^3(c)$ with the Gaussian curvature $K=c$. 
\end{rem}

\begin{rem}{\textbf{(Assumptions)}}
Throughout this paper, we are going to assume the followings:
\begin{itemize}
\item  The space time $L^4_1(f,c)$ has no open part with constant sectional curvature, i.e., 
$$\frac{f''(t)}{f(t)} -\frac{f'(t)^2 + c}{f(t)^2}\neq 0$$
 for all $t\in I$.

\item $T$ does not vanish on $M$. 

\item The relative null bundle $\mathcal{N}$ is a smooth 1-dimensional distribution.

\item $M$ is  connected and it has no open part lying on a 3-dimensional totally geodesic hypersurface of $L^4_1(f,c)$.
\end{itemize} 
\end{rem}

On the other hand,   
by considering Lemma \ref{LemmaLn1f0CurvTensor}, one can observe that the Codazzi and Ricci equations turn into 
\begin{eqnarray}    
    \label{CodEqMostGeneral} \left(\langle X,X\rangle Y_0-\langle X,Y\rangle X_0\right)\left(-\frac{f''}{f} + \frac{(f')^2 + c}{f^2}\right)\eta &=& \left(\nabla^\perp_{X} h\right)(Y, X) - \left(\nabla^\perp_{Y} h\right)(X, X),
\end{eqnarray}
where  $\nabla^\perp h$ is defined by
$$\left(\nabla^\perp_X h\right)(Y, Z) = \nabla^\perp_X h(Y, Z) - h(\nabla_X Y, Z) - h(Y, \nabla_X Z).$$
Moreover, the Ricci equation turns into
\begin{eqnarray}
    \label{RicciLn1fc} 
    R^\perp(X,Y)\xi &=& h(X, A_\xi Y) - h(A_\xi X, Y).
\end{eqnarray}
Consequently, $M$ has flat normal bundle, i.e., $R^\perp$ if and only if all shape operators of $M$ have the same casual type, \cite{ChenJVanderVeken2007}.


\section{Surfaces with Positive Relative Nullity}

\subsection{Space-like Surfaces}
Let $M$ be a space-like surface in $L^4_1(f,c)$. In this case, by considering the decomposition \eqref{RWTS-etaandTDef}, we define an orthonormal frame field $\{e_1,e_2;e_3,e_4\}$ and an `angle'  function $\theta$ by
\begin{equation}\label{RWTS-etaandTDefNewTimel}
\left.\frac{\partial}{\partial t} \right|_M=\sinh\theta \, e_1+\cosh\theta \, e_3.  
\end{equation}

\begin{proposition}\label{PropPRNL41fcSpacelike}
Let $M$ be a space-like surface in $L^4_1(f,c)$. Then, $M$ has positive relative nullity if and only if  its second fundamental form satisfies 
\begin{equation}\label{PRN2ndFundForm}
h(e_1,e_1)=h(e_1,e_2)=0,\qquad h(e_2,e_2)=\xi,
\end{equation}
where $e_1,e_2$ are orthonormal tangent vector fields defined by \eqref{RWTS-etaandTDefNewTimel} and $\xi$ is a normal vector field.
\end{proposition}

\begin{proof}
In order to prove the necessary condition, assume that $M$ has positive relative nullity. Then, there exists an orthonormal base $\{\hat e_1,\hat e_2\}$ of the tangent bundle of $M$ such that 
\begin{equation}\label{PropPRNL41fcSpacelikeEq1}
h(\hat e_1,\hat e_1)=h(\hat e_1,\hat e_2)=0,\qquad h(\hat e_2,\hat e_2)=\xi
\end{equation}
for a normal vector field $\xi$. By using the Codazzi equation \eqref{CodEqMostGeneral} for $X=\hat e_1$ and $Y=\hat e_2$, we get
\begin{equation}\label{PropPRNL41fcSpacelikeEq2}
\left(-\frac{f''}{f} + \frac{(f')^2 + c}{f^2}\right)\beta\cosh\theta e_3= \langle \nabla_{\hat e_1}\hat e_1,\hat e_2\rangle \xi,
\end{equation}
where we put $\beta=(\hat e_2)_0$. Now, consider the open subset $$\mathcal O=\{p\in M|\beta(p)\neq0\}$$ of $M$. We are going to prove that $\mathcal O$ is empty.
Because of \eqref{PropPRNL41fcSpacelikeEq2}, we have $\xi=\alpha e_3$ on $\mathcal O$ from which, along with \eqref{AhRelatedBy}, we get
\begin{equation}\label{PropPRNL41fcSpacelikeEq3}
\left.A_{e_4}\right|_{\mathcal O}=0,
\end{equation}
where $\alpha$ is a smooth function.

On the other hand, from \eqref{RWTS-etaandTDefNewTimel} we get
\begin{equation}\label{PropPRNL41fcSpacelikeEq4}
\langle\widetilde\nabla_X\frac{\partial}{\partial t},e_4 \rangle=\langle\widetilde\nabla_X\left(\sinh\theta \, e_1+\cosh\theta \, e_3\right),e_4\rangle
\end{equation}
whenever $X$ is a tangent vector field. By a further computation considering Lemma \ref{LemmaLn1f0LCConnect}, \eqref{RWTS-etaandTDefNewTimel} and \eqref{PropPRNL41fcSpacelikeEq4}, we obtain
\begin{equation}\label{PropPRNL41fcSpacelikeEq5}
\langle\left(\sinh\theta \, h(e_1,X)+\cosh\theta \, \nabla^\perp_X e_3\right),e_4\rangle=0
\end{equation}
from which we get $ \nabla^\perp e_3=0$ because of \eqref{PropPRNL41fcSpacelikeEq3}. Since the codimension of $M$ is 2, we also have$ \nabla^\perp e_4=0$. Therefore, \eqref{PropPRNL41fcSpacelikeEq3} implies
$$\widetilde\nabla e_4=0 \mbox{ on $\mathcal O$}.$$
Therefore, if $\mathcal O$ is non-empty, then any component of $\mathcal O$ lies on a totally geodesic hypersurface of $L^4_1(f,c)$. 
By the hypothesis, this is a contradiction. Hence, $\mathcal O$ is empty which yields that $\hat e_1=\pm e_1$. This completes the proof of the necessary condition.

Converse is obvious because  \eqref{PRN2ndFundForm} implies $(e_1)_p\in\mathcal{N}_p$ for all $p\in M$.
\end{proof}


\begin{lemma}\label{PRNSpacelikeLemma1}
Let $M$ be a space-like surface in $L^4_1(f,c)$. 
Then, $M$ has positive relative nullity if and only if the vector fields $e_1,e_2;e_3,e_4$ and the function $\theta$ defined by \eqref{RWTS-etaandTDefNewTimel} satisfy

\begin{subequations}\label{PRNSpacelikeCond}
\begin{eqnarray}
\label{PRNSpacelikeCondEq1} \widetilde\nabla_{e_1}e_1=0, &\qquad& \widetilde\nabla_{e_2}e_1=\omega e_2,\\
\label{PRNSpacelikeCondEq2} \widetilde\nabla_{e_1}e_2=0, &\qquad& \widetilde\nabla_{e_2}e_2=-\omega e_1-h^3_{22}e_3+h^4_{22}e_4,\\
\label{PRNSpacelikeCondEq3} \widetilde\nabla_{e_1}e_3=0, &\qquad& \widetilde\nabla_{e_2}e_3=-h^3_{22}e_2,\\
\label{PRNSpacelikeCondEq4} \widetilde\nabla_{e_1}e_4=0, &\qquad& \widetilde\nabla_{e_2}e_4=-h^4_{22}e_2, \\
\label{PRNSpacelikeCondEq5} e_1(\theta)=\frac{f'}f\cosh\theta,&\qquad&e_2(\theta)=0
\end{eqnarray}
\end{subequations}
for some smooth functions $\omega,\ h^3_{22}$ and $h^4_{22}$.
\end{lemma}
\begin{proof}
In order to prove the necessary condition, assume that $M$ has positive relative nullity. Then, by Proposition \ref{PropPRNL41fcSpacelike}, the equation \eqref{PRN2ndFundForm} is satisfied. Then, we have
\eqref{RWTS-etaandTDefNewTimel} implies
\begin{equation}\label{PRNSpacelikeLemma1Eq0}
A_{e_3}=\left(\begin{array}{cc}0&0\\0&h^3_{22}
\end{array}\right),\qquad 
A_{e_4}=\left(\begin{array}{cc}0&0\\0&h^4_{22}
\end{array}\right)
\end{equation}
because of \eqref{AhRelatedBy}, where we put $h^\alpha_{22}=\langle\xi,e_\alpha\rangle,\ \alpha=3,4$.
Note that 
\eqref{RWTS-etaandTDefNewTimel} implies
\begin{align}\label{PRNSpacelikeLemma1Eq1}
\begin{split}
\widetilde\nabla_{e_1}\frac{\partial}{\partial t}=&\widetilde\nabla_{e_1}\left(\sinh\theta \, e_1+\cosh\theta \, e_3\right),\\
\widetilde\nabla_{e_2}\frac{\partial}{\partial t}=&\widetilde\nabla_{e_2}\left(\sinh\theta \, e_1+\cosh\theta \, e_3\right).
\end{split}
\end{align}
By using Lemma \ref{LemmaLn1f0LCConnect}, \eqref{PRN2ndFundForm} and \eqref{PRNSpacelikeLemma1Eq0}, we get
\begin{align}\label{PRNSpacelikeLemma1Eq2}
\begin{split}
\frac{f'}f(\cosh^2\theta e_1&+\sinh\theta \cosh\theta e_3)=e_1(\theta)\zeta+\sinh\theta \nabla_{e_1}e_1+\cosh\theta\nabla^\perp_{e_1}e_3,\\
\frac{f'}fe_2&=e_2(\theta)\zeta+\sinh\theta \nabla_{e_2}e_1+\cosh\theta(-h^3_{22}e_2+\nabla^\perp_{e_2}e_3),
\end{split}
\end{align}
where $\zeta=\cosh\theta \, e_1+\sinh\theta \, e_3$. 
By a further computation using \eqref{PRNSpacelikeLemma1Eq0} and \eqref{PRNSpacelikeLemma1Eq2}, we obtain \eqref{PRNSpacelikeCond}, for the function $\omega$ defined by 
$$\omega=\langle \nabla_{e_2}e_1,e_2\rangle.$$
This completes the proof of the necessary condition. Converse is obvious.
\end{proof}


\subsection{Time-like Surfaces} In this subsection, we are going to consider time-like surfaces with positive relative nullity. First, we obtain the analogous of Proposition \ref{PropPRNL41fcSpacelike} for time-like surfaces.
\begin{proposition}\label{PRNProp1}
Let $M$ be a time-like surface in $L^4_1(f,c)$. Then, $M$ has positive relative nullity if and only if its second fundamental form satisfies \eqref{PRN2ndFundForm}, where $e_1,e_2$ form an orthonormal base for the tangent bundle of $M$ such that $e_1$ is proportional to $ \left(\frac{\partial}{\partial t}\right)^T$ and $\xi$ is  a normal vector field.
\end{proposition}

\begin{proof}
In order to prove the necessary condition, assume that $M$ has positive relative nullity. Note that
if $\eta=0$ on $M$, then $e_1=\frac{\partial}{\partial t}$ and $e_2$ satisfies
$$\widetilde\nabla_{e_1}e_1=0\mbox{ and }\widetilde\nabla_{e_2}e_1=\frac{f'}{f}e_2$$
because of Lemma \ref{LemmaLn1f0LCConnect} (See also \cite[Lemma 3.3]{ChenJVanderVeken2007}). Therefore, the proof is completed for this case.

Now, we assume the existence of $p\in M$ such that  $\eta(p)\neq 0$ and consider the open subset $\mathcal O=\{q\in M|  \eta(q)\neq0\}$.  In this case, the vector fields $T$ and $\eta$ defined by \eqref{RWTS-etaandTDef} are time-like and space-like, respectively. First, we are going to prove that $\mathcal N_p$ can not be degenerated. Towards contradiction, assume that $\mathcal{N}_p$ is degenerated on $\mathcal O$. 
In this case, there exists a base field $\{X,Y\}$ of the tangent bundle of $\mathcal O$ such that
\begin{equation*}
h(X,X)=h(X,Y)=0,\quad h(Y,Y)=\xi,\qquad \langle X,X\rangle=\langle Y,Y\rangle=0,\quad \langle X,Y\rangle=-1.
\end{equation*}
Consequently, the Codazzi equation \eqref{CodEqMostGeneral} implies
\begin{eqnarray*}    
    X_0\left(-\frac{f''}{f} + \frac{(f')^2 + c}{f^2}\right)\eta &=& 0
\end{eqnarray*}
from which we obtain $X_0=\langle X,T\rangle=0$. However, this is a contradiction because $X$ is light-like and $T$ is time-like.  

Consequently, one can define a local orthonormal frame field $\{e_1,e_2;e_3,e_4\}$ on $\mathcal O$ and a  function $\theta$ by
\begin{equation}\label{RWTS-TimeLikeetaandTDefNewLocal}
\left.\frac{\partial}{\partial t} \right|_{\mathcal O}=\cosh\theta \, e_1+\sinh\theta \, e_3.  
\end{equation}
Similar to the proof of Proposition \eqref{PropPRNL41fcSpacelike}, we assume that there exists an orthonormal base $\{\hat e_1,\hat e_2\}$ of the tangent bundle of $\mathcal O$ such that \eqref{PropPRNL41fcSpacelikeEq1} is satisfied and $\langle\hat e_1,\hat e_1\rangle=\varepsilon\in\{-1,1\}$. By an analogous computation, we obtain $\hat e_1=\pm e_1$ to complete the proof of the necessary condition. Converse is obvious.
\end{proof}


\begin{lemma}\label{PRNTimelikeLemma0}
Let $M$ be a time-like surface in $L^4_1(f,c)$. Then, there exists a constant $a$ such that
\begin{equation}\label{PRNTimelikeCond0}
\theta=\sinh^{-1}\frac{a}f
\end{equation}
on $M$.
\end{lemma}

\begin{proof}
If  $\eta=0$ on $M$, then \eqref{PRNTimelikeCond0} is satisfied for $a=0$. Assume the existence of $p\in M$ such that $\eta(p)\neq0$. Similar to the proof of Proposition \ref{PRNProp1}, we consider the connected component $\mathcal O$ of $\{q\in M|  \eta(q)\neq0\}$ such that $p\in \mathcal O$. Let  $\{e_1,e_2;e_3,e_4\}$ be the local orthonormal frame field  on $\mathcal O$ defined by \eqref{RWTS-TimeLikeetaandTDefNewLocal}.

Similar to the proof of Lemma \ref{PRNSpacelikeLemma1}, we consider the equation
\begin{align*}
\begin{split}
\widetilde\nabla_{e_1}\frac{\partial}{\partial t}=&\widetilde\nabla_{e_1}\left(\cosh\theta \, e_1+\sinh\theta \, e_3\right),\\
\widetilde\nabla_{e_2}\frac{\partial}{\partial t}=&\widetilde\nabla_{e_2}\left(\cosh\theta \, e_1+\sinh\theta \, e_3\right),
\end{split}
\end{align*}
 Lemma \ref{LemmaLn1f0LCConnect}, \eqref{PRN2ndFundForm} and \eqref{PRNSpacelikeLemma1Eq0}, to get
\begin{eqnarray}
\label{PRNTimelikeLemma0Eq1} \nabla_{e_1}e_1=0, &\qquad& \nabla_{e_2}e_1=\omega e_2,\\
\label{PRNTimelikeLemma0Eq2} e_1(\theta)=-\frac{f'}f\sinh\theta,&\qquad& e_2(\theta)=0.
\end{eqnarray}
on $\mathcal O$, where $\omega$ is a smooth function.

Note that \eqref{PRNTimelikeLemma0Eq1} and \eqref{PRNTimelikeLemma0Eq2} imply
$[\frac1{\cosh\theta} e_1,G e_2]=0$
if $G$ is a smooth function satisfying 
\begin{equation}\label{PRNTimelikeLemma0Eqq2}
e_1(G)=G\omega.
\end{equation}
Therefore, for any $q\in \mathcal O$ there exists a local coordinate system $(\mathcal O_q,(u,v))$ exists such that $\mathcal O_q\ni q$ and 
\begin{equation}\label{PRNTimelikeLemma0Eqq1}
e_1=\cosh\theta\partial_u\qquad e_2=\frac 1G\partial_v.
\end{equation}
Consequently, \eqref{PRNTimelikeLemma0Eq2} gives
$$\theta=\theta(u),\qquad \theta'(u)\coth\theta(u)+\frac{f'(u)}{f(u)}=0$$
on $\mathcal O_q$. Therefore, we have 
$$\theta(u)=\sinh^{-1}\frac{a}{f(u)}$$
for a constant $a$ on $\mathcal O_p$. Since $M$ is connected, the continuity of $\theta$ implies $\mathcal O=M$ and \eqref{PRNTimelikeCond0} on $M$.
%
%
%
\end{proof}


Consider the case $a=0$, i.e., $\eta=0$ on $M$.  Note that the classification of surfaces satisfying this condition is obtained in \cite{ChenJVanderVeken2007}:
\begin{lemma}\label{ChenJVanderVekenLemma3.3}\cite{ChenJVanderVeken2007}. 
 If $M$ is a surface of $L^4_1(f,c)$ such that $\frac{\partial}{\partial t} $ is everywhere tangent to $M$, then
$M$ is an open part of $I\times_f\alpha\subset I\times_fR^3(c)$, where $\alpha=\alpha(s)$ is a curve in $R^3(c)$. Moreover, if $\alpha(s)$  is of unit
speed, then the second fundamental form of $M$ in $L^4_1(f,c)$  satisfies
$$h(\partial_t,\partial_t)=h(\partial_t,\partial_s)=0,\qquad h(\partial_s,\partial_s)=\alpha''(s),$$
where we have identified the vector $\alpha''(s):=\nabla^{R^3(c)}_{\alpha'(s)}\alpha'$ in $R^3(c)$ with its vertical lift in $L^4_1(f,c)$  via $\Pi^2$.
\end{lemma}
As a direct consequence of this lemma, we have the following corollary.
\begin{corry}
If $M$ is a time-like surface in $L^4_1(f,c)$ with $\eta=0$, then $M$ has positive relative nullity. 
\end{corry}


Next, we are going to consider the case $a\neq0$. In this case, by considering the decomposition \eqref{RWTS-etaandTDef}, we define an orthonormal frame field $\{e_1,e_2;e_3,e_4\}$ and an `angle'  function $\theta$ by
\begin{equation}\label{RWTS-TimeLikeetaandTDefNew}
\left.\frac{\partial}{\partial t} \right|_M=\cosh\theta \, e_1+\sinh\theta \, e_3.  
\end{equation}

The proof of the following lemma is similar to the proof of Lemma \ref{PRNSpacelikeLemma1}.
\begin{lemma}\label{PRNTimelikeLemma1}
Let $M$ be a time-like surface in $L^4_1(f,c)$ with $\eta\neq0$. 
Then, $M$ has positive relative nullity if and only if the vector fields $e_1,e_2;e_3,e_4$ 
and the function $\theta$ defined by \eqref{RWTS-TimeLikeetaandTDefNew} satisfy
\begin{subequations}\label{PRNTimelikeCond}
\begin{eqnarray}
\label{PRNTimelikeCondEq1} \widetilde\nabla_{e_1}e_1=0, &\qquad& \widetilde\nabla_{e_2}e_1=\omega e_2,\\
\label{PRNTimelikeCondEq2} \widetilde\nabla_{e_1}e_2=0, &\qquad& \widetilde\nabla_{e_2}e_2=\omega e_1+h^3_{22}e_3+h^4_{22}e_4,\\
\label{PRNTimelikeCondEq3} \widetilde\nabla_{e_1}e_3=0, &\qquad& \widetilde\nabla_{e_2}e_3=-h^3_{22}e_2,\\
\label{PRNTimelikeCondEq4} \widetilde\nabla_{e_1}e_4=0, &\qquad& \widetilde\nabla_{e_2}e_4=-h^4_{22}e_2
\end{eqnarray}
\end{subequations}
for some smooth functions $\omega,\ h^3_{22}$ and $h^4_{22}$.
\end{lemma}


\section{Surfaces in $L^4_1(f,0)$}
In this section, we are going to consider non-degenerated surfaces in $L^4_1(f,0)$.

First, we construct a space-like surface with positive relative nullity.
\begin{proposition}\label{SpSurfacePRNExample}
Let  $a_1,a_2,\phi_1: I_v \to\mathbb R$ be some smooth functions and assume the existence of an $a\in\mathbb R$ such that $a^2-f(u)^2>0$ whenever $u\in I_u$,  where $I_u\subset I$ and $I_v$ are open intervals. Let $\phi_2,\phi_3$ and $\alpha_1,\alpha_2,\alpha_3:I_v\to\mathbb R^3$ some smooth functions satisfying
\begin{equation}\label{SpSurfacePRNExampleEq1}
\phi_2'=-a_1\phi_1,\qquad \phi_3'=-a_2\phi_1,
\end{equation}
and
\begin{equation}\label{SpSurfacePRNExampleEq2}
\left\{
\begin{array}{rcll}
\alpha_1'&=&a_1(v)\alpha_2+a_2(v)\alpha_3,\qquad & \alpha_1(v_0)=C_1,\\
\alpha_2'&=&-a_1(v)\alpha_1,& \alpha_2(v_0)=C_2,\\
\alpha_3'&=&-a_2(v)\alpha_1,& \alpha_3(v_0)=C_3,
\end{array}
\right.
\end{equation}
where $v_0\in I_v$ and $C_1,C_2,C_3$ are constant orthonormal vectors in $\mathbb E^3$. 

Then, the space-like surface $M$ in $L^4_1(f,0)$ parametrized by
\begin{align}\label{SpSurfacePRNExampleEq3}
\begin{split}
\phi(u,v)=&\left(u,\phi_1(v)\alpha_1(v)+\left(\int_{u_0}^u\frac{a}{f(\bar u)\sqrt{a^2-f(\bar u)^2}}d\bar u+\phi_2(v)\right)\alpha_2(v)
+\phi_3(v)\alpha_3(v)\right)\\
\end{split}
\end{align}
has positive relative nullity.
\end{proposition}
\begin{proof}
By a direct computation, we observe that the vector fields $e_1,e_2,e_3,e_4$ given by
\begin{eqnarray*}
e_1=\frac{f}{\sqrt{a^2-f^2}}\partial_u, &\qquad& e_2=\frac{1}{\sqrt G}\partial_v,\\
e_3=\left(\frac{a}{f(u)},\frac{\sqrt{a^2-f(u)^2}}{f(u)}\alpha _2\right),&\qquad& e_4=\left(0,\frac1f \alpha_3\right)
\end{eqnarray*}
form the orthonormal frame field defined by \eqref{RWTS-etaandTDefNewTimel}, where we put
$$G=f^2 \left(a_1 \left(\int_{u_0}^u \frac{a}{f(\bar u) \sqrt{a^2-f(\bar u)^2}} \, d\bar u+{\phi _2}\right)+a_2 {\phi _3}-{\phi _1}'\right)^2.$$
By a further computation, we obtain that \eqref{PRNSpacelikeCond} is satisfied for
\begin{align}\label{SpSurfacePRNExampleEqCond}
\begin{split}
\omega=&\frac{f' \sqrt{G \left(a^2-f^2\right)}+a a_1 f}{\sqrt Gf^2},\\
h^3_{22}=&\frac{a_1 f\sqrt{a^2-f^2}+a \sqrt{G}  f'}{\sqrt{G} f^2},\\
h^4_{22}=&\frac{a_2}{\sqrt{G}}.
\end{split}
\end{align}
Hence, $M$ has positive relative nullity because of Lemma \ref{PRNSpacelikeLemma1}.
\end{proof}


Next, we prove the following classification theorem:
\begin{theorem}\label{ThmL41f0ProofEq1}
Let $M$ be a space-like surface in $L^4_1(f,0)$. Then, $M$ has positive relative nullity if and only if it is locally congruent to the surface described in Proposition \ref{SpSurfacePRNExample}.
\end{theorem}

\begin{proof}
The sufficient condition is proved in Proposition \ref{SpSurfacePRNExample}. In order to prove the necessary condition, we assume that $M$ is a space-like surface in $L^4_1(f,0)$ with positive relative nullity,  $p\in M$ and let $\phi$ denote the position vector of $M$. Also, we put $\tilde\phi=\Pi^2\circ\phi$. Note that Lemma \ref{PRNSpacelikeLemma1} implies that \eqref{PRNSpacelikeCond} is satisfied. From \eqref{PRNSpacelikeCondEq1}, \eqref{PRNSpacelikeCondEq2} and \eqref{PRNSpacelikeCondEq5} we have
$$[-\frac{1}{\sinh\theta}e_1,Ge_2]=0$$
for a function $G$ satisfying $e_1(G)-G\omega=0$. Therefore, there exists a local coordinate system $(\mathcal O_p,(u,v))$ such that
\begin{subequations}\label{ThmL41f0ProofEq1ALL}
\begin{eqnarray}
\label{ThmL41f0ProofEq1a}e_1=-\sinh\theta\partial_u&=&-\sinh\theta(1,\tilde\phi_u),\\
\label{ThmL41f0ProofEq1b}e_2=\frac 1G\partial_v&=&\frac 1G(0,\tilde\phi_v),
\end{eqnarray}
\end{subequations}
where \eqref{ThmL41f0ProofEq1a} follows from \eqref{RWTS-etaandTDefNewTimel}. Note that \eqref{ThmL41f0ProofEq1ALL} implies
\begin{equation}\label{ThmL41f0ProofEq1yer}
\phi(u,v)=(u,\tilde\phi(u,v))
\end{equation}
on $\mathcal O_p$.

On the other hand, because of \eqref{PRNSpacelikeCondEq5}, on $\mathcal O_p$ we have
$$\theta=\theta(u),\qquad \theta'(u)\sinh\theta(u)+\frac{f'(u)}{f(u)}\cosh\theta(u)=0$$
from which we get
$$\cosh\theta(u)=\frac{a}{f(u)},\qquad \sinh\theta(u)=\frac{\sqrt{a^2-f(u)^2}}{f(u)}.$$
Consequently, \eqref{ThmL41f0ProofEq1a} turns into
\begin{equation}
\label{ThmL41f0ProofEq2a}e_1=-\frac{\sqrt{a^2-f(u)^2}}{f(u)}(1,\tilde\phi_u).
\end{equation}
Next, by considering \eqref{ThmL41f0ProofEq2a}, we use Lemma \ref{LemmaLn1f0LCConnect} to get
\begin{equation}
\label{ThmL41f0ProofEq3}\widetilde\nabla_{e_1}e_1=\left(0,\frac{\left(a^2-2 f ^2\right) f'  \tilde\phi_u+f  \left(a^2-f ^2\right) \tilde\phi_{uu}}{f ^3}\right).
\end{equation}
By combining \eqref{ThmL41f0ProofEq3} with the first equation of \eqref{PRNSpacelikeCondEq1}, we obtain
\begin{equation}
\label{ThmL41f0ProofEq4} \tilde\phi_u(u,v)=\frac{a}{f(u) \sqrt{a^2-f(u)^2}}\alpha_2(v)
\end{equation}
for a smooth $\mathbb E^3$-valued function $\alpha_2$.

On the other hand, we use Lemma \ref{LemmaLn1f0LCConnect} and \eqref{ThmL41f0ProofEq1a} to observe that the second equation in \eqref{PRNSpacelikeCondEq1} implies
\begin{equation}
\label{ThmL41f0ProofEq5} \tilde\phi_v(u,v)=\frac{G(u,v)}{f(u)}\alpha_1(v)
\end{equation}
for a smooth $\mathbb E^3$-valued function $\alpha_1$. Since $\{e_1,e_2\}$ is orthonormal, \eqref{ThmL41f0ProofEq1ALL} implies that $\{\alpha_1,\alpha_2\}$ is orthonormal. 

Now, let $\alpha_3$ be a  $\mathbb E^3$-valued function  so that
$$g_c(\alpha_i,\alpha_j)=\delta_{ij},\qquad i,j=1,2,3.$$
Then, we have
\begin{equation}
\label{ThmL41f0ProofEq6} \tilde\phi(u,v)=\phi_1(u,v)\alpha_1(v)+\hat\phi_2(u,v)\alpha_2(v)+\phi_3(u,v)\alpha_3(v)
\end{equation}
and
\begin{equation}\label{ThmL41f0ProofEq7}
\left\{
\begin{array}{rcl}
\alpha_1'&=&a_1(v)\alpha_2+a_2(v)\alpha_3\\
\alpha_2'&=&-a_1(v)\alpha_1+a_3(v)\alpha_3,\\
\alpha_3'&=&-a_2(v)\alpha_1-a_3(v)\alpha_2,
\end{array}
\right.
\end{equation}
for some smooth functions $\phi_1,\hat\phi_2,\phi_3,a_1,a_2,a_3$. By a further computation using \eqref{ThmL41f0ProofEq4}-\eqref{ThmL41f0ProofEq7}, we obtain $a_3=0$,  \eqref{SpSurfacePRNExampleEq1} and
\begin{equation}\label{ThmL41f0ProofEq8}
\partial_u\phi_1=\partial_u\phi_3=0,\qquad \hat\phi_2(u,v)=\int_{u_0}^u\frac{a}{f(\bar u)\sqrt{a^2-f(\bar u)^2}}d\bar u+\phi_2(v).
\end{equation}
By combining \eqref{ThmL41f0ProofEq6} and \eqref{ThmL41f0ProofEq8} with \eqref{ThmL41f0ProofEq1yer}, we get \eqref{SpSurfacePRNExampleEq3}. Moreover, because of $a_3=0$,  \eqref{ThmL41f0ProofEq7} turns into \eqref{SpSurfacePRNExampleEq2}. Hence, $M$ is congruent to the surface described in Proposition \ref{SpSurfacePRNExample}.
\end{proof}


In the remaining part, we deal with time-like surfaces in $L^4_1(f,0)$. 
Next, we obtain analogous of Proposition \ref{SpSurfacePRNExample} for time-like surfaces.
\begin{proposition}\label{TiSurfacePRNExample}
Let  $a_1,a_2,\phi_1:I_v\to\mathbb R$ be some smooth functions and assume that $\phi_2,\phi_3$ and $\alpha_1,\alpha_2,\alpha_3:I_v\to\mathbb R^3$ are smooth functions satisfying \eqref{SpSurfacePRNExampleEq1} and \eqref{SpSurfacePRNExampleEq2} for a $v_0\in I_v$ and some constant orthonormal vectors $C_1,C_2,C_3$ in $\mathbb E^3$.

Then, the time-like surface $M$ in $L^4_1(f,0)$ parametrized by
\begin{equation}\label{TiSurfacePRNExampleEq3}
\phi(u,v)=\left(u,\phi_1(v)\alpha_1(v)+\left(\int_{u_0}^u\frac{a}{f(\bar u)\sqrt{a^2+f(\bar u)^2}}d\bar u+\phi_2(v)\right)\alpha_2(v)
+\phi_3(v)\alpha_3(v)\right)
\end{equation}
has positive relative nullity, where $a$ is a non-zero constant.
\end{proposition}

\begin{proof}
By a direct computation, we observe that $M$ is time-like and \eqref{PRNTimelikeCond0} is satisfied for $a\neq0$. Similar to the proof of Proposition \ref{SpSurfacePRNExample}, we obtain that the vector fields $e_1,e_2,e_3,e_4$ defined by \eqref{RWTS-TimeLikeetaandTDefNew} satisfy \eqref{PRNTimelikeCond} for 
\begin{align}\label{TiSurfacePRNExampleEqCond}
\begin{split}
\omega=&\frac{f' \sqrt{G \left(a^2+f^2\right)}+a a_1 f}{\sqrt G f^2},\\
h^3_{22}=&\frac{a_1 f \sqrt{a^2+f^2}+a \sqrt{G} f'}{\sqrt G f^2},\\
h^4_{22}=&\frac{a_2}{\sqrt G},
\end{split}
\end{align}
where we put
$$G=f^2 \left(a_1 \left(\int_{{u0}}^u \frac{a}{f(z) \sqrt{a^2+f(z)^2}} \, dz+{\phi _2}\right)+a_2 {\phi _3}-{\phi _1}'\right)^2.$$
Hence, $M$ has positive relative nullity because of Lemma \ref{PRNTimelikeLemma1}. 
\end{proof}

Next, we obtain the following classification theorem for time-like surfaces in $L^4_1(f,0)$.  
\begin{theorem}
Let $M$ be a time-like surface in $L^4_1(f,0)$. Then, $M$ has positive relative nullity if and only if it is one of the following surfaces:
\begin{enumerate}
\item[(1)] The surface described in Lemma \ref{ChenJVanderVekenLemma3.3} for a curve $\alpha$ in $\mathbb E^3$,
\item[(2)] The surface locally congruent to the surface described in Proposition \ref{TiSurfacePRNExample}.
\end{enumerate}
\end{theorem}

\begin{proof}
Similar to the proof of Theorem \ref{ThmL41f0ProofEq1}, assume that $M$ is a time-like surface in $L^4_1(f,0)$ with positive relative nullity,  $p\in M$ and let $\phi$ denote the position vector of $M$. By Lemma \ref{PRNTimelikeLemma0}, \eqref{PRNTimelikeCond0} is satisfied for a constant $a$. If $a=0$, then we have $\eta=0$ on $M$ and Lemma \ref{ChenJVanderVekenLemma3.3} implies that $M$ is the surface  given in  the Case (1) of the theorem. 

Now, we consider the case $a\neq0$. Then, Lemma  \ref{PRNTimelikeLemma1} implies that \eqref{PRNTimelikeCond} is satisfied. Let $(\mathcal O_p,(u,v))$ the coordinate system constructed in the proof of Lemma \ref{PRNTimelikeLemma0} satisfying \eqref{PRNTimelikeLemma0Eqq1} for a function $G$ such that \eqref{PRNTimelikeLemma0Eqq2}. Moreover, by Lemma \ref{PRNTimelikeLemma1}, \eqref{PRNTimelikeCond} is satisfied. Similar to the proof of Theorem \ref{ThmL41f0ProofEq1}, by considering the equations in \eqref{PRNTimelikeCondEq1}, we use Lemma \ref{LemmaLn1f0LCConnect} to get
\begin{equation}
\label{ThmL41f0ProofEq1Ti} \tilde\phi_u(u,v)=\frac{a}{f(u) \sqrt{a^2+f(u)^2}}\alpha_2(v)
\end{equation}
and \eqref{ThmL41f0ProofEq5} for some smooth $\mathbb E^3$-valued functions $\alpha_1$ and $\alpha_2$. Consequently, we get \eqref{ThmL41f0ProofEq6} on $\mathcal O_p$. By a further computation considering \eqref{ThmL41f0ProofEq5} and \eqref{ThmL41f0ProofEq1Ti},  we obtain that $M$ is the surface given in Case (2) of the theorem. This completes the proof of the necessary condition.

The converse is proved in Lemma \ref{ChenJVanderVekenLemma3.3}  and Proposition \ref{TiSurfacePRNExample}.
\end{proof}


\section{Surfaces in Lorentzian Product Spaces}
In this section, we are going to consider the case $f=1$ and $c=\pm 1$ to study surfaces in the product spaces $\mathbb E^1_1\times R^3(c)$.

Let $M$ be a non-degenerated surface in the Cartesian product space $L^4_1(1,c)=\mathbb E^1_1\times R^3(c)$. 
In this section $\hat\nabla$ stands for the Levi-Civita connection of the flat ambient space $\mathbb E^5_r$, 
where $r = 1$ if $c = 1$ and $r = 2$ if $c = -1$. Note that $\hat\nabla$ and $\widetilde\nabla$  are related with
\begin{subequations}\label{LCConnectionsRelatedAll}
\begin{eqnarray}
\label{LCConnectionsRelated1}\hat\nabla_X Y&=&\widetilde\nabla_X Y-c\left( \langle X,Y \rangle + \langle X,T\rangle\langle Y,T\rangle\right)e_5\\
\label{LCConnectionsRelated2}\hat\nabla_X \xi&=&\widetilde\nabla_X \xi-c\langle X,T\rangle\langle \xi,\eta\rangle e_5
\end{eqnarray}
whenever $X,Y$ are tangent to $M$ and $\xi$ is a normal vector field, where $e_5$ denotes the restriction of the unit normal vector field of $\mathbb E^1_1\times R^{3}(c)$ in  $\mathbb E^{5}_r$ into $M$. Furthermore, we have
\begin{equation}\label{LCConnectionsRelated3}
\hat\nabla_X e_5=(X+\langle X,T\rangle T)+\langle X,T\rangle \eta.
\end{equation}
\end{subequations}

Now, assume that $M$ has positive relative nullity. In this case, since $f=1$, \eqref{PRNSpacelikeCondEq5} in Lemma \ref{PRNSpacelikeLemma1} and Lemma \ref{PRNTimelikeLemma0} imply
\begin{equation}\label{SectPRNE11R3c}
\theta=\theta_0 
\end{equation}
for a constant $\theta_0$, where $\theta$ is defined by \eqref{RWTS-etaandTDefNewTimel} (resp. \eqref{RWTS-TimeLikeetaandTDefNew}) if $M$ is space-like (resp. time-like). Note that we have $\theta_0 \neq0$ when $M$ is space-like.

\subsection{Surfaces in $\mathbb E^1_1\times\mathbb S^3$}
In this subsection, we construct two families of surfaces in $\mathbb E^1_1\times\mathbb S^3$ with positive relative nullity.
\begin{proposition}\label{SpPRNinE11S3Example}
Let  $a_1,a_2,a_3:I_v\to\mathbb R$ be some smooth functions and $\alpha_1,\alpha_2,\alpha_3,\alpha_4:I_v\to\mathbb R^4$ are smooth functions satisfying
\begin{equation}\label{SpPRNinE11S3ExampleEq1}
\left\{
\begin{array}{rcll}
\alpha_1'&=&a_1(v)\alpha_3,& \alpha_1(v_0)=C_1,\\
\alpha_2'&=&a_2(v)\alpha_3,& \alpha_2(v_0)=C_2,\\
\alpha_3'&=&-a_1(v)\alpha_1-a_2(v)\alpha_2+a_3(v)\alpha_4,& \alpha_3(v_0)=C_3,\\
\alpha_4'&=&-a_3(v)\alpha_3,& \alpha_4(v_0)=C_4,\\
\end{array}\right.
\end{equation}
where $v_0\in I_v$ and $C_1,C_2,C_3,C_4$   are constant, orthonormal vectors in $\mathbb E^4$.

Then, the space-like surface $M$ in $ \mathbb E^1_1\times\mathbb S^3$ parametrized by
\begin{equation}\label{SpPRNinE11S3ExampleEq2}
\phi(u,v)=\left(u\sinh \theta_0,\cos(u\cosh \theta_0)\alpha_1(v)+\sin(u\cosh \theta_0)\alpha_2(v) \right)
\end{equation}
has positive relative nullity, where  $\theta_0$ is a non-zero constant.
\end{proposition}

\begin{proof}
 By a direct computation, we obtain that \eqref{PRNSpacelikeCond} is satisfied for  
\begin{align}\label{SpPRNinE11S3ExampleEq3}
\begin{split}
\omega=&\frac{\cosh \theta_0 (a_1 \sin (u \cosh \theta_0 )-a_2 \cos (u \cosh \theta_0 ))}{a_1 \cos (u \cosh \theta_0 )+a_2 \sin (u \cosh \theta_0 )},\\
h^3_{22}=&\frac{\sinh \theta_0 (a_2 \cos (u \cosh \theta_0 )-a_1 \sin (u \cosh \theta_0 ))}{a_1 \cos (u \cosh \theta_0 )+a_2 \sin (u \cosh \theta_0 )},\\
h^4_{22}=&-\frac{a_3}{a_1 \cos (u \cosh \theta_0 )+a_2 \sin (u \cosh \theta_0 )}.
\end{split}
\end{align}
Hence, $M$ has positive relative nullity because of Lemma \ref{PRNSpacelikeLemma1}.
\end{proof}


\begin{proposition}\label{TiPRNinE11S3Example}
Let  $a_1,a_2,a_3:I_v\to\mathbb R$ be some smooth functions and $\alpha_1,\alpha_2,\alpha_3,\alpha_4:I_v\to\mathbb R^4$ are smooth functions satisfying
\eqref{SpPRNinE11S3ExampleEq1} for a $v_0\in I_v$ and some constant, orthonormal vectors $C_1,C_2,C_3,C_4$  in $\mathbb E^4_1$.

Then, the time-like surface $M$ in $ \mathbb E^1_1\times\mathbb S^3$ parametrized by
\begin{equation}\label{TiPRNinE11S3ExampleEq2}
\phi(u,v)=\left(u\cosh \theta_0,\cos(u\sinh \theta_0)\alpha_1(v)+\sin(u\sinh \theta_0)\alpha_2(v) \right)
\end{equation}
has positive relative nullity, where  $\theta_0$ is a constant.
\end{proposition}

\begin{proof}
If $\theta_0=0$, then \eqref{TiPRNinE11S3ExampleEq2} turns into
$$\phi(u,v)=\left(u,\alpha_1(v)\right).$$ 
In this case, by Lemma \ref{ChenJVanderVekenLemma3.3}, $M$ has positive relative nullity. 

Assume that $\theta_0\neq0$. By a direct computation, we obtain that \eqref{PRNTimelikeCond} is satisfied for 
\begin{align}\label{TiPRNinE11S3ExampleEq3}
\begin{split}
\omega=&\frac{\sinh \theta_0 (a_2 \cos (u \sinh \theta_0)-a_1 \sin (u \sinh \theta_0))}{a_1 \cos (u \sinh \theta_0)+a_2 \sin (u \sinh \theta_0)},\\
h^3_{22}=&\frac{\cosh \theta_0 (a_2 \cos (u \sinh \theta_0)-a_1 \sin (u \sinh \theta_0))}{a_1 \cos (u \sinh \theta_0)+a_2 \sin (u \sinh \theta_0)},\\
h^4_{22}=&-\frac{a_3}{a_1 \cos (u \sinh \theta_0)+a_2 \sin (u \sinh \theta_0)}.
\end{split}
\end{align}
Hence, $M$ has positive relative nullity because of Lemma \ref{PRNTimelikeLemma1}.
\end{proof}


\subsection{Surfaces in $\mathbb E^1_1\times\mathbb H^3$}
In this subsection, we construct two families of surfaces in $\mathbb E^1_1\times\mathbb H^3$ with positive relative nullity.
\begin{proposition}\label{SpPRNinE11H3Example}
Let  $a_1,a_2,a_3:I_v\to\mathbb R$ be some smooth functions and $\alpha_1,\alpha_2,\alpha_3,\alpha_4:I_v\to\mathbb R^4$ are smooth functions satisfying
\begin{equation}\label{SpPRNinE11H3ExampleEq1}
\left\{
\begin{array}{rcll}
\alpha_1'&=&a_1(v)\alpha_3,& \alpha_1(v_0)=C_1,\\
\alpha_2'&=&a_2(v)\alpha_3,& \alpha_2(v_0)=C_2,\\
\alpha_3'&=&a_1(v)\alpha_1-a_2(v)\alpha_2+a_3(v)\alpha_4,& \alpha_3(v_0)=C_3,\\
\alpha_4'&=&-a_3(v)\alpha_3,& \alpha_4(v_0)=C_4,
\end{array}\right.
\end{equation}
where $v_0\in I_v$, $C_1$ is timelike, $C_2,C_3,C_4$ are spacelike constant, orthonormal vectors in $\mathbb E^4_1$.

Then, the space-like surface $M$ in $ \mathbb E^1_1\times\mathbb H^3$ parametrized by
\begin{equation}\label{SpPRNinE11H3ExampleEq2}
\phi(u,v)=\left(u\sinh \theta_0,\cosh(u\cosh \theta_0)\alpha_1(v)+\sinh(u\cosh \theta_0)\alpha_2(v) \right)
\end{equation}
has positive relative nullity, where  $\theta_0$ is a non-zero constant.
\end{proposition}

\begin{proof}
The surface $M$ satisfies \eqref{PRNSpacelikeCond} for 
\begin{align}\label{SpPRNinE11H3ExampleEq3}
\begin{split}
\omega=&-\frac{\cosh \theta_0 (a_1 \sinh (u \cosh \theta_0 )+a_2 \cosh (u \cosh \theta_0 ))}
{a_1 \cosh (u \cosh \theta_0 )+a_2 \sinh (u \cosh \theta_0 )},\\
h^3_{22}=&\frac{\sinh \theta_0 (a_1 \sinh (u \cosh \theta_0 )+a_2 \cosh (u \cosh \theta_0 ))}
{a_1 \cosh (u \cosh \theta_0 )+a_2 \sinh (u \cosh \theta_0 )},\\
h^4_{22}=&-\frac{a_3}{a_1 \cosh (u \cosh \theta_0 )+a_2 \sinh (u \cosh \theta_0 )}.
\end{split}
\end{align}
Hence, $M$ has positive relative nullity because of Lemma \ref{PRNSpacelikeLemma1}.
\end{proof}


\begin{proposition}\label{TiPRNinE11H3Example}
Let  $a_1,a_2,a_3:I_v\to\mathbb R$ be some smooth functions and $\alpha_1,\alpha_2,\alpha_3,\alpha_4:I_v\to\mathbb R^4$ 
are smooth functions satisfying
\eqref{SpPRNinE11H3ExampleEq1} for a $v_0\in I_v$ and some constant, orthonormal vectors $C_1,C_2,C_3,C_4$  in $\mathbb E^4_1$.
Then, the time-like surface $M$ in $ \mathbb E^1_1\times\mathbb H^3$ parametrized by
\begin{equation}\label{TiPRNinE11H3ExampleEq2}
\phi(u,v)=\left(u\cosh \theta_0,\cosh(u\sinh \theta_0)\alpha_1(v)+\sinh(u\sinh \theta_0)\alpha_2(v) \right)
\end{equation}
has positive relative nullity, where  $\theta_0$ is a non-zero constant.
\end{proposition}

\begin{proof}
The surface $M$ satisfies \eqref{PRNTimelikeCond} for 
\begin{align}\label{TiPRNinE11H3ExampleEq3}
\begin{split}
\omega=&\frac{\sinh \theta_0 (a_1 \sinh (u \sinh \theta_0)+a_2 \cosh (u \sinh \theta_0))}{a_1 \cosh (u \sinh \theta_0)+a_2 \sinh (u \sinh \theta_0)},\\
h^3_{22}=&\frac{\cosh \theta_0 (a_1 \sinh (u \sinh \theta_0)+a_2 \cosh (u \sinh \theta_0))}{a_1 \cosh (u \sinh \theta_0)+a_2 \sinh (u \sinh \theta_0)},\\
h^4_{22}=&-\frac{a_3}{a_1 \cosh (u \sinh \theta_0)+a_2 \sinh (u \sinh \theta_0)}.
\end{split}
\end{align}
Hence, $M$ has positive relative nullity because of Lemma \ref{PRNSpacelikeLemma1}.
\end{proof}

\begin{rem} From the proofs of Proposition \ref{SpPRNinE11S3Example}-Proposition \ref{TiPRNinE11H3Example}, 
it can be easily seen that $h^4_{22}=0$ for $a_3=0$. In this case, the surfaces has open part lying on 
a 3-dimensional totally geodesic hypersurface of $L^4_1(f,c)$. Thus, the interrior of
$$\{v_0\in I_v| a_3(v_0)=0\}$$
must be empty.
\end{rem}

\subsection{Local Classification Theorem}
In this subsection, we obtain the following local classification of surfaces in $L^4_1(1,c)$ with positive relative nullity.
\begin{theorem}\label{THME11S31}
Let $M$ be a non-degenerated surface in $L^4_1(1,c)$ with positive relative nullity. Then, the followings hold:
\begin{itemize}
\item[(i)]  If $M$ is a space-like surface and $c=1$, then it is locally congruent to the surface described in Proposition \ref{SpPRNinE11S3Example}.
\item[(ii)] If $M$ is a time-like surface and $c=1$, then it is locally congruent to the surface described in Proposition \ref{TiPRNinE11S3Example}. 
\item[(iii)]  If $M$ is a space-like surface and $c=-1$, then it is locally congruent to the surface described in Proposition \ref{SpPRNinE11H3Example}.
\item[(iv)] If $M$ is a time-like surface and $c=-1$, then it is locally congruent to the surface described in Proposition \ref{TiPRNinE11H3Example}. 
\end{itemize}
\end{theorem}

\begin{proof}
First, we are going to prove (i) and (iii) of the theorem by considering a space-like surface $M$ in $\mathbb E^1_1\times\mathbb S^3$ 
with positive relative nullity. Let $\phi$ be the position vector of $M$. Then by combining  the first equations in \eqref{PRNSpacelikeCondEq1}, \eqref{PRNSpacelikeCondEq2}and \eqref{PRNSpacelikeCondEq4} with \eqref{LCConnectionsRelatedAll} we get
\begin{subequations}\label{THME11S31Eq1ALL}
\begin{eqnarray}
\label{THME11S31Eq1a} \hat\nabla_{e_1}e_1&=&-c\left(\cosh^2{\theta_0}\right)e_5, \\
\label{THME11S31Eq1b} \hat\nabla_{e_1}e_2&=&0, \\
\label{THME11S31Eq1c} \hat\nabla_{e_1}e_4&=&0. 
\end{eqnarray}
\end{subequations}

On the other hand, by combining \eqref{SectPRNE11R3c} with \eqref{RWTS-etaandTDefNewTimel}, we obtain
\begin{equation}\label{THME11S31Eq4}
\left.\frac{\partial}{\partial t} \right|_M=\sinh\theta_0 \, e_1+\cosh\theta_0 \, e_3.  
\end{equation}

Now, let $p\in M$. Similar to the proof of Theorem \ref{ThmL41f0ProofEq1}, by considering $[e_1,Ge_2]=(e_1(G)-G\omega)e_2$, we choose a local coordinate system $(\mathcal O_p,(u,v))$ such that
\begin{subequations}\label{THME11S31Eq5ALL}
\begin{eqnarray}
\label{THME11S31Eq5a}e_1=\partial_u&=&(-\sinh\theta_0,\tilde\phi_u),\\
\label{THME11S31Eq5b}e_2=\frac 1G\partial_v&=&\frac 1G(0,\tilde\phi_v),
\end{eqnarray}
\end{subequations}
where $G$ is a smooth function satisfying $e_1(G)-G\omega=0$ and $p\in \mathcal O_p$. Note that because of \eqref{THME11S31Eq5ALL}, up to a suitable translation on $u$, we have 
\begin{equation}\label{THME11S31Eq6}
\phi(u,v)=(-u\sinh\theta_0,\tilde\phi(u,v))
\end{equation}
on $\mathcal O_p$, where $\tilde\phi=\Pi^2(\Phi)$. 

Next, by a direct computation considering \eqref{THME11S31Eq1ALL}, \eqref{THME11S31Eq5ALL} and \eqref{THME11S31Eq6}, we get
\begin{subequations}\label{THME11S31Eq7ALL}
\begin{eqnarray}
\label{THME11S31Eq7a}\tilde\phi_{uu}+c\cosh^2\theta_0 \tilde\phi&=&0,\\
\label{THME11S31Eq7b}e_2&=&(0,\alpha_3(v)),\\
\label{THME11S31Eq7c}e_4&=&(0,\alpha_4(v)),\\
\end{eqnarray}
\end{subequations}
where $\alpha_3,\alpha_4:I_v\to\mathbb R^4$ are some smooth functions and $I_v$ is an open interval. By solving \eqref{THME11S31Eq7a}, we get
\begin{equation}\label{THME11S31Eq8c=1}
\tilde\phi(u,v)=\cos(u\sinh\theta)\alpha_1(v)+\sin(u\sinh\theta)\alpha_2(v)
\end{equation}
for $c=1$ and
\begin{equation}\label{THME11S31Eq8c=-1}
\tilde\phi(u,v)=\cosh(u\sinh\theta)\alpha_1(v)+\sinh(u\sinh\theta)\alpha_2(v)
\end{equation}
for $c=-1$, where $\alpha_1,\alpha_2:I_v\to\mathbb R^4$ are some smooth functions. By combining \eqref{THME11S31Eq6} with \eqref{THME11S31Eq8c=1} and \eqref{THME11S31Eq8c=-1}, we get \eqref{SpPRNinE11S3ExampleEq2} and \eqref{SpPRNinE11H3ExampleEq2}, respectively. Moreover, since $\{e_1,e_2,e_4\}$ is an orthonormal set, $\{\alpha_1,\alpha_2,\alpha_3,\alpha_4\}$ must be orthonormal in $\mathbb E^4_{r-1}$. By a further computation we obtain \eqref{SpPRNinE11S3ExampleEq1} (resp. \eqref{SpPRNinE11H3ExampleEq1}) for $c=1$ (resp. $c=-1$). Hence we have (i) and (iii) of the theorem. (ii) and (iv) is proved by a similar technique. 
\end{proof}

\section*{Acknowledgements}
This work was carried out during the 1001 project supported by the Scientific and Technological Research Council of T\"urkiye (T\"UB\.ITAK)  (Project Number: 121F352).

\end{document}